\documentclass[12pt]{article}
\usepackage{amssymb}
\usepackage{amsmath,amsfonts,amsthm}
\usepackage{srcltx}
\begin{document}

                             \begin{center}
                             	
                                     \bf{ COMMENTS ON "ORBITS OF AUTOMOPHISM GROUPS OF FIELDS"}\\
                                            Pramod K. Sharma\\
                                             email: pksharma1944@yahoo.com\\
                                             Department of Mathematics, Sikkim Central University, \\                                             
                                             Gangtok-73710, INDIA                                     
                                     \end{center} 
          \begin{center}                             
          \bf {ABSTRACT}     \end{center} Let $R$ be a commutative $k-$algebra over a field $k$. Assume $R$ is a noetherian, infinite, integral domain. The group of $k-$automorphisms of $R$,i.e.$Aut_k(R)$ acts in a natural way on $(R-k)$.In the first part of this article, we study the structure of $R$ when the orbit space $(R-k)/Aut_k(R)$ is finite.We note that most of the results, not particularly relevent to fields, in [1,\S 2] hold in this case as well. Moreover, we prove that $R$ is a field. In the second part, we study a special case of the Conjecture 2.1 in [1] :  If $K/k$ is a non trivial field extension where $k$ is algebraically closed and $\mid (K-k)/Aut_k(K) \mid = 1$ then $K$ is algebraically closed. In the end, we give an elementary proof of [1,Theorem 1.1] in case $K$ is finitely generated over its prime subfield.\\
          
           \section{Introduction}
          
          Let $ K/k$ be a non trivial field extension. Authors in $ [1, \S2]$ conjecture that the orbit space $ (K-k)/Aut_k(K) $ is finite if and only if either both $K$ and $k$ are finite or both are algebraically closed. From the results of the authors, it is clear that $K$ is finite if and only if $k$ is finite. Moreover,if $K$ is algebraically closed then so is $k$. The converse is open, and  several results are proved in $[1]$ for this case. In this note, we prove that if $R$ is an infinite, noetherian integral domain which is  an algebra over a field $k$ such that $\mid (R-k)/Aut_k(R)\mid < \infty$ then most  of the results, not particularly relevent to fields, in [1,\S 2] hold in this case as well. We, infact, note that $R$ is a field in case charateristic of $R$ is $ p > 0$ or $R$ is integrally closed. Further,we note that if $K$ is algebraically closed then $\mid(K-k)/Aut_k(K)\mid = 1$. Note that if the conjecture is true, $\mid(K-k)/Aut_k(K)\mid = 1$ should imlpy that $ K$ is algebraically closed. We are not able to prove this. However, we observe that that a field $E$ is algebraically closed if and only if $ f(E) = E$ for all $ f(X) \in E[X]$, and prove that $ f(K) = K$ for a class polynomials including all non zero polynomials in $ k[X] $. It is also proved that if $\mid(K-k)/Aut_k(K)\mid < \infty$, then for any non constant $f(X) \in K[X]$ with degree $ \leq 2, f(K) = K$.\\
          
          \section{Main Results} 
          	
          	Throughout,we assume that $R$ is an infinite commutative $k-$algebra over a field $k$ which is a noetherian integral domain such that $\mid (R-k)/Aut_k(R)\mid < \infty.$ \\
          	
          	Theorem 2.1. The field $k$ is infinite and integrally closed in $R$.\\
          	
          	Proof. If $\mid k \mid < \infty$, then as $\mid (R-k)/Aut_k(R)\mid < \infty$, $ \mid (R/Aut_k(R)\mid < \infty.$  Hence by $[2,
          	\mbox { Corollary } 16]$, $R$ is a finite field.This contradicts our assumption that $R$ is infinite. Hence $ \mid k \mid = \infty$. Now, let $ \alpha \in (R-k)$ be integral over $k$. Then for each $a\in k, a\alpha $ is integral over $k$, moreover $ \{ a\alpha \mid  a\in k\}$ is an infinite subset of $(R-k)$. Note that if  $ \beta \in (R-k)$ is integral over $k$, then for any $ \sigma \in Aut_k(R), \sigma (\beta)$ is integral over $k$. Hence orbit of $\beta$, i.e. $O(\beta) = \{\sigma(\beta) \mid \sigma \in Aut_k(R \}$ is a finite set. This implies  $\mid (R-k)/Aut_k(R)\mid = \infty$, a contradiction to the assumption that  $\mid (R-k)/Aut_k(R)\mid < \infty.$ Hence $k$ is infinite and integrally closed in $R$. \\
          	           	
          	Theorem 2.2. If characteristic of $k$ is $p> 0$, then $ k^p =k$, and $R^p = R$. Moreover, $R$ is a field.\\
          	
          	Proof. By   $ \mbox{ Theorem 2.1 }$, $k$ is intergrally closed in $R$. Hence $ (R-k)^p  \subset (R-k)$. Consequently $$  (R-k) \supset  (R-k)^p \supset \ldots  \supset  (R-k)^{p^m} \supset \ldots  $$ is a chain of orbit closed subsets of  $(R-k)$ under the action of $Aut_k(R)$. Since  $\mid (R-k)/Aut_k(R)\mid <\infty $, there exists $ n \geq 1$ such that $$  (R-k)^{p^n}=  (R-k)^{p^{(n+1)}}.$$ Thus for any $ \lambda \in  (R-k)$, there exists $\mu \in  (R-k)$ such that  
          	
          	                                  \begin{eqnarray*}
          	                                       \lambda^{p^n} & = & \mu^{p^{(n+1)}}\\
          	                                       \Rightarrow (\lambda - \mu^p)^{p^n} & = & 0\\
          	                                       \Rightarrow  \lambda & = & \mu^p \\
          	                                       \Rightarrow (R-k) & = &  (R-k)^p \subset R^p          	                                                                
          	                                  \end{eqnarray*}
          Now, if $\lambda \in (R-k), a\in k$, then $$ \lambda, (\lambda - a) \in (R-k) = (R-k)^p \subset R^p.$$ Assume $ \lambda = \alpha^p,  \lambda - a = \beta^p \mbox { for } \alpha, \beta \in R$. Then \begin{eqnarray*} 
          (\alpha - \beta)^p & = & \alpha^p - \beta^p = a \in R^p\\
         \Rightarrow & k & \subset R^p \\
          \Rightarrow  R & = & R^P \mbox { since }  (R-k) = (R-k)^p.                    
          \end{eqnarray*}  	 Now as $R =R^P$ and 
          $k$ is integrally closed, $ k = k^p$. We shall now show that $R$ is a field. Let $ I$ be a non-zero radical ideal in $R$. As $R = R^p$, for any $ a\in I$, there exists $ b \in R$ such that $ a= b^p$. As $I$ is a radical ideal $ b \in I$. Consequently $ I = I^p$, i.e. $I^{p-1}I = I$. As $R$ is noetherian, there exists $ \alpha \in I^{p-1}$ such that $ (1-\alpha)I = (0)$. This is not possible as $ I \neq (0)$ as well as $ \alpha \neq 1$. Hence $R$ has no non-zero radical ideal. Thus $R$ is a field,.\\
          
          Theorem 2.3. For any $ x \in (R-k)$ and $ c\in k$, there exists $ \sigma \in  Aut_k(R)$ such that $ \sigma(x) = x + c$. Hence $ k + x \subset  O(x)$. \\
          
          Proof. The proof is similar to the proof of [1, lemma 2.6]. \\
          
         Theorem 2.4. If $R$ is integrally closed, then $ (R-k)^l = (R-k) $ for all $ l \geq 1$. Moreover, $R$ is a field, and $R^l = R$ as well as $ k^l = k$.\\
         
         Proof. It suffices to prove the statement assuming $l$ is prime. In view of Theorem 2.2, we can assume that $l$ is other than the characteristic of $k$ if that is prime. As in theorem 2.2, $$  (R-k) \supset  (R-k)^l \supset \ldots  \supset  (R-k)^{l^m} \supset \ldots  $$ is a chain of orbit closed subsets of $(R-k)$. As $\mid (R-k)/Aut_k(R)\mid < \infty$, there exists $ n \geq 1$ such that $$ (R-k)^{l^m } = (R-k)^{l^{(m+1)}} $$                      
         for all $ m \geq n$. Thus for any $\lambda \in (R-k)$, there exists $ \mu \in (R-k)$ such that \begin{eqnarray*} 
         	  \lambda^{l^m} & = & \mu^{l^{(m+1)}}\\
         	\Rightarrow  (\lambda \mu^{-l})^{l^m} & = & 1 \\
         	\Rightarrow    \lambda \mu^{-l} &  \in R                     
         \end{eqnarray*}	                      
         since $R$ is integrally closed. Further, as $k$ is integrally closed in $R, \lambda \mu^{-l} \in (k^*)_{l^m}$, the subgroup of $(l^m)th $ roots of unity in $k$. Therefore $ \lambda \in (R-k)^l(k^*)_{l^m}$ for all $ m \geq n$. Consequemtly $(R-k) \subset (R-k)^l(k^*)_{l^m}$. Next, for any $ m \geq n$  $$ (R-k)^{l^m} =  (R-k)^{l^{2m}}.$$ Hence as above, we conclude that for any $\lambda \in (R-k), \lambda \in (R-k)^{l^m} (k^*)_{l^m}$. Thus $(R-k) \subset (R-k)^{l^m}(k^*)_{l^m}$. We, now, consider two cases:\\
         Case 1. $(k^*)_{l^m}  \subsetneqq (k^*)_{l^{(m+1)}}.$ \\
         
         In this case, for any $c\in (k^*)_{l^m}$, there exists $b \in (k^*)_{l^{(m+1)}}$ such that $ c = b^l$. Hence,as $ (R-k) \subset (R-k)^l(k^*)_{l^m}$, we have $$ (R-k) \subset (R-k)^l(k^*)_{l^m} =  (R-k)^l \subset (R-k).$$  Consequently $ (R-k) =(R-k)^l$.\\
         
         Case 2. $(k^*)_{l^m} =(k^*)_{l^{(m+1)}}$ for all $ m \geq n.$\\
         
         We have $$  (R-k) \subset (R-k)^l(k^*)_{l^m} \subset (R-k)k^* \subset (R-k).$$ Consequently $ (R-k) = (R-k)^l(k^*)_{l^m}$ for all  $ m \geq n.$ Further, since $$ (R-k) \subset (R-k)^{l^m}(k^*)_{l^m} \subset (R-k)k^* \subset (R-k)$$  \begin{eqnarray*}
         	\Rightarrow  (R-k) & = & (R-k)^{l^m}(k^*)_{l^m} \\  \Rightarrow  (R-k) & = & (R-k)^{l^m}(k^*)_{l^{(m+1)}} \mbox{ since } (k^*)_{l^m} =(k^*)_{l^{(m+1)}} \\ \Rightarrow (R-k)^l & = & (R-k)^{l^{(m+1)}}(k^*)_{l^{(m+1)}}
         	\end{eqnarray*} If $ m > n$, then $ m-1 \geq n$. Hence $$(R-k)^l = (R-k)^{l^m}(k^*)_{l^m} = (R-k).$$ We shall now show that $R$ is a field. Let $ I$ be a non-zero radical ideal in $R$, then as $ I\cap k = (0), I-(0) \subset (R-k).$ Let $ a(\neq 0)\in I,$ then $ a = \lambda^l$ for some $ \lambda \in (R - k)$  since $R-k)^l = (R-k)$. Thus as $I$ is radical ideal, $ \lambda \in I.$ Consequently  
         	 \begin{eqnarray*} I^l & = & I\\
         	 \Rightarrow  I^2 & = & I \\ 
         	 \Rightarrow  (1-\lambda) I & = & (0) \end{eqnarray*} for some $\lambda \in I$. As $ \lambda \neq 0 ,1$, it is a non trivial idempotent . This is not possible, hence $I = (0)$. Therefore $R$ is a field. The last part of the statement follows by [ 1, Proposition 2.10]\\ 
         	
         	Remark 2.5. If $U$ is the group of units in $R$, and $ U \cap (R-k) \neq  \phi$, then $ R^l = R$ and $ k^l = k$. This gives an alternative proof of the last part of the statement.\\
         	
         	We shall first prove that $ U^l = U$. Let $ v \in U \cap (R-k) =  U \cap (R-k)^l$. Then $ v = \lambda ^l$ for some $\lambda \in (R-k)$. Clearly $ \lambda \in U$. Hence $ v \in U^l$. Therefore $ U \cap (R-k) \subset U^l $. Next, note that  $$ U = (U\cap (R-k) \cup (U\cap k^*) = ( U\cap (R-k)) \cup k^* \subset U^l \cup k^*$$  Thus to prove $ U = U^l$, it suffices to show that $ k^*  \subset U^l$. Let $ a\in k^*$ and $ \lambda \in  U \cap (R-k)$. Then $ \lambda a \in U \cap (R-K) .$ Thus $ \lambda^{-1} (\lambda a)= a \in U^l $ since $ ( U \cap (R-k))\subset U^l$. Hence $ k^* \subset U^l$. Thus $ U = U^l$. Now, let $ a \in k^*$. Then there exists $ b \in U$ such that $ b^l = a$. As $k$ is integrally closed in $R, b \in k^*. $  Hence $ k = k^l$. This implies
         	  \begin{eqnarray*}
         	      R^l & = & (R-k)^l \cup k^l\\
         	      & = & (R-k) \cup k\\
         	      & = & R
         	  \end{eqnarray*}  Hence the assertion holds.\\

             Theorem 2.6. If $R$ is integrally closed, then $ k = R^{Aut_k(R)} = \{\lambda \mid \sigma (\lambda) = \lambda \mbox { for all }\sigma \in Aut_k(R)\}$.\\  
             
             Proof. Let $ \lambda \in R^{Aut_k(R)}, \lambda \notin k.$ By [2, Lemma 5], $ \lambda $ is a unit. Therefore $ L = R^{Aut_k(R)}$ is a field containing $k$. By Theorem 2.4, $(R-k)^l = (R-k)$ for any $ l \geq 1$. Thus since $ \lambda \notin k, X^l - \lambda$ has no root in $k$. Moreover, by Theorem 2.3, for any $ a\in k$ and a root $\mu$ of $X^l - \lambda, \mu + a $ is also a root of $ X^l - \lambda$ since for any $\sigma \in R^{Aut_k(R)}, \sigma(\mu)$ is also a root of $ X^l -\lambda$. As $ \mid k \mid = \infty$, this is not possible. Hence $ k = R^{Aut_k(R)}$.\\
             
             Remark 2.7. (i) If the charateristic of $k$ is $ p > 0$, then we can drop the condition that $R$ is integrally closed.This can be seen by taking $ l= p$.   \\
             (ii) Under the conditions of the theorem, $\mid O(\lambda)\mid <  \infty$ if and only if $ \lambda \in k$. If $ O(\lambda) = \{\lambda_1, \ldots , \lambda_t \}$, then $ p(X) = (X-\lambda_1)\ldots(X-\lambda_t)\in k[X]$. Thus each $ \lambda_i, i=1,\ldots t$ is integral over $k$, and consequently $ \lambda \in k$. The converse is clear.\\
             
             Theorem 2.8. If $\lambda \in (R-k)$, then $ S_{\lambda} = \{a \in k^* \mid \sigma(\lambda)  = a\lambda \mbox { for some } \sigma \in Aut_k(R)\}$ is a subgroup of finite index in $k^*$.\\
             
             Proof. Let $a,b \in S_{\lambda}$. Then there exist $ \sigma, \tau \in Aut_k(R)$ such that $ \sigma(\lambda) = a\lambda, \tau(\lambda) = b\lambda$. Therefore $ \sigma\tau(\lambda)=ab \lambda$ and $\sigma^{-1}(\lambda) = a^{-1}\lambda $. Hence $ ab, a^{-1} \in S_{\lambda}$. Thus  $S_{\lambda}$ is a subgroup of $ k^*$. Assume $ [k^* :  S_{\lambda}]= \infty$. Choose an infinite set $ \{b_1, \ldots, b_n, \ldots\}$ in $k^*$ such that $ b_i S_{\lambda} \neq b_j S_{\lambda}$ for all $ i \neq j$. We claim $ O(b_i\lambda) \neq  O(b_j\lambda)$ whenever $ i\neq j$. If not, then there exist $ i \neq j$ such that   \begin{eqnarray*}  O(b_i\lambda) & = & O(b_j\lambda) \\ \Rightarrow \sigma(b_i\lambda ) & = & b_j\lambda \mbox{ for some } \sigma \in Aut_k(R)\\ \Rightarrow \sigma(\lambda) & = & b_i^{-1}b_j \lambda \\ \Rightarrow  b_i^{-1}b_j & \in & S_{\lambda} \\ \Rightarrow b_i \lambda & = & b_j \lambda  \end{eqnarray*} As   $ b_i S_{\lambda} \neq b_j S_{\lambda}$ for all $ i \neq j$ , the claim follows. This cotradicts the assumption that $ \mid (R-k)/ Aut_k(R) \mid < \infty $.  Hence  $ [k^* :  S_{\lambda}] < \infty$.\\
             
             Remark 2.9. (i) Let $k$ be algebraically closed, and let  $ [k^* :  S_{\lambda}] =m < \infty$.  Then since $ (k^*)^m = k^*$, $S_{\lambda} = k^*$. Hence for any $ c\in k^*$, there exists $ \sigma \in Aut_k(R)$ such that $ \sigma(\lambda) = c \lambda$.\\
             (ii) Assume $R$ is integrally closed and 
             $ (R-k) \cap U \neq \phi$, where $U$ is the group of units of $R$. Then also the assertion of the theorem holds, i.e. we need not assume that $k$ is algebraically closed in this case. This follows since $ (k^*)^m = k^*$ by remark 2.5.\\  
             
             \bfseries{Theorem 2.10.}  Let $k$ be algebraically closed and $ a \neq(0), b \in k$. Then $$ L_{a,b} :( R - k) \longrightarrow (R - k)$$$$ x \longrightarrow  ax + b $$ is one-one, onto map. Moreover $ G = \{ L_{a, b} \mid a\neq (0), b \in k\} $ is a group with respect to composition of maps. Further this group acts on the orbits of $ (R - k)$ under the action of $ Aut_k(R) $ trivially.\\
             
             Proof. One can easily check that $G$ is a group and since the action of $ L_{a,b}$ commutes with the action of $ Aut_k(R)$, the group $G$ acts on the orbits of $ (R - k)$. Note that $ L_{a,b}(x) = L_{1, b}. L_{a,0}(x)$ for any $ x \in (R - k)$, thus $ L_{a,b} = L_{1,b} L_{a,0}$. Now, by the Theorem 2.3 and the Remark 2.9(i), it is clear that $L_{a,0} $ and $ L_{1,b}$ act trivially on the orbits $ (R - k)$, the result follows.\\    
             
             \section{A look at Conjecture} 
            
                       We shall first give an elementary proof of [1, Theorem 1.1] in case  the field $K$ is finitely generated over its prime subfield. Then we study a particular case of the Cojecture 2.1 in [1]:\\
                       "Let $K/k$ be a non trivial extension of fields. Then the number of orbits of $Aut_k(K)$  on $(K-k)$ is finite if and only if either both $K \mbox{ and } k $ are finite or both are algebrically closed."\\
                       The authors in [1] have proved that if the number of orbits of $Aut_k(K)$  on $(K-k)$ is finite, then $K$ is finite if and only if $k$ is finite.Further, it is noted that if $K$ is algebraically closed then so is $k$. Thus it remains to show that under the given condition if $k$ is algebraically closed then so is $K$.
                       Based on the conjecture we ask : If $K/k$ is a non trivial field extension where $k$ is algebraically  closed, then is it true that K is algebraically closed if and only if $Aut_k(K)$  on $(K-k)$ has one orbit? Before we look into this, we prove: \\                     
             	       
             	       Theorem 3.1. Let $K$ be a field with $\mid K/Aut(K)\mid < \infty$. If $K$ is finitely generated over its prime subfield, then $K$ is finite.\\
             	       
             	       Proof. It is noted in [1]that characteristic of $K$ is $ p> 0$ and it is perfect i.e., $ K^p = K$. Let $\mathbb{F}_p$ be the prime subfield of $K$. As $K$ is finitely generated over $\mathbb{F}_p$, $K$ has finite transcendence degree over $\mathbb{F}_p$. Let $S$ be a transcdence basis of $K \mid \mathbb{F}_p$. Then $K \mid \mathbb{F}_p (S)$ is finite algebraic. If $S = \phi$, then clearly $K$ is finite. Further if $ S \neq \phi$, then as $K$ is perfect $ K \neq \mathbb{F}_p(S)$. As $K$ is perfect, the Fr$\ddot{o}$benius endomorphism of $K$ $\sigma$(say) is an automorphism. Therefore as $ [K: \mathbb{F}_p(S)] < \infty $, $ \sigma ( \mathbb{F}_p (S) = \mathbb{F}_p(S)$. This however is not true. Consequently $ S = \phi$, and $K$ is finite. \\
             	       
             	       Here after, we assume that $K/k$ is a non trivial field extension where $k$ is algebraically closed.\\
             	       
             	       Lemma 3.2. If $K$ is algebraically closed then  action of $Aut_k(K)$ over $(K-k)$ has one orbit.\\
             	       
             	       Proof. Let $ x,y \in (K-k)$. Since $k$ is algebraically closed, $x,y$ are transcendental over $k$. Choose a transcendental basis $S$ of $K/k$ containing x, and a transcendental basis $T$ of $K/k$ containing y. As $S$ and $T$ have same cardinality, there exists a bijection from $S$ to $T$ mapping x to y. This extends to a $k-$ isomorphism $\sigma$ from the field $k(S)$ to the field $k(T)$ which maps x to y. As $K$ is algebraic closure of $k(S)$ as well as $k(T)$, $\sigma$ extends to an automorphism $\tau$ of $K$ such that $ \tau(x) = y$. Hence the result follows.\\  
             	       
             	       Lemma 3.3. A field $K $ is algebraically closed if and only if for every non-constant polynomial $f(X) \in K[X], f(K) = \{f(\lambda) \mid \lambda \in K \} = K$.\\
             	       
             	       Proof. Let $K$ be algebraicallty closed  and $f(X) \in K[X]$ be a non-constant polynomial. Note that for any $ \lambda \in K, g(X) = f(X) - \lambda$ is a non constant polynomial in $K[X]$. Hence has a root in $K$. If $a$ is a root of $g(X)$ in $K$, then $f(a) = \lambda$. Therefore $ f(K) = K$. Conversely, let $ f(X) \in K[X]$ be any non-constant polynomial. Then since $ f(K) = K$, there exists $ b \in K$ such that $ f(b) = 0$. Hence $K$ is algebraically closed.\\
             	       
             	       Remark 3.4. (i) For a field $K$ to be algebraically closed, it is sufficient to assume that $p(K) = K$ for every irreducible polynomial $p(X) \in K[X]$.\\  
             	       
             	       (ii) For a polynomial $f(X) \in K[X]$, the condition $ f(K) = K$ is not equivalent to the fact that $f(X)$ splits over $K$, e.g. if $ K = \mathbb{R}$(the field of reals), and $ f(X) = X^3 -1$, then $ f(\mathbb{R}) = \mathbb{R}$, but $f(X)$ does not splir over $ \mathbb{R}$. Moreover if $ f(X) = X^2 - 1$, then $ f(\mathbb{R}) \neq \mathbb{R} $, but this splits.\\
             	       
             	       (iii) A field $E$ is algebraically closed if and only if every prime ideal in $ E[X]$ is of the form $ (X - \alpha)E[X] (\alpha \in E)$.\\

             	       Theorem 3.5. If $Aut_k(K)$ has one orbit over $(K-k)$, then for any non-constant $ f(X) \in k[X]$, $ f(K) = K$, moreover, $f(K-k)= K-k$.\\
             	       
             	       Proof. First of all, note that since $k$ is algebraically closed  and $ f(X)$ is non-constant, for any $ a\in k$ , there exists $ b \in k $ such that $ f(b) = a.$ Now, let $\alpha \in (K-k)$ be any element, then  $f(\alpha)= \beta \in (K-k)$. Thus, since  $Aut_k(K)$ has one orbit over $(K-k)$, for any $z \in (K-k)$, there exists $ \sigma \in Aut_k(K)$ such that $ \sigma (\beta) = z$ . Hence $ \sigma (f(\alpha))= f(\sigma (\alpha) = z$. Consequently $ f(K-k) = (K-k)$. Therefore $ f(K) = K$. Finally, for any $ \alpha \in (K-k), f(\alpha)\in (K-k)$ since $ \alpha$ is transcendental over $k$ and $f(X) \in k[X]$. Hence as $ f(K) = K, f(K-k) = K-k$.   \\
             	       
             	       Remark 3.6. (i) If $ f(X) \in k[X]$, and $ \lambda (\neq 0) \in K$, then $(\lambda f(X)) (K) = K$ as well as $ (\lambda + f(X))(K) = K$.\\
             	       
             	       (ii) If for $ f(X), g(X) \in K[X], f(K) = K = g(K)$, then for $ h(X) = f(X) o g(X) =  f(g(X)), h(K) = K. $ Moreover for any $ \lambda (\neq(0)) \in K, \lambda f(K) = K  \mbox {and also } (f(X) + \lambda)(K) = K$. \\ 
             	       
             	       (iii) By the Theorem 3.5, $ K^l = K \forall l \geqq 1$. Thus $K$ and $k$ are perfect. Moreover, if for $f(X) \in K[X], f(K) = K$, then $ f^l(K) = K$ and also $ f(K^l) = K.$\\
             	       
             	       (iv) If for $ f(X) \in K[X], f(K) = K$, then for any $ p(X) \in k[X], p(f)(K) = K .$ Thus for any $ a\neq (0), b \in K, \mbox { and } h(x) = p(aX + b), h(K)= K.$
             	       \\
             	       
             	       (v) If for $ f(X) \in K[X], f(K) = K$, then for any $ \sigma \in  Aut_k(K), f^{\sigma}(X)(K) = K$ where $ f^{\sigma}(X)$ denotes the polynomial obtained from $ f(X)$ by applying $ \sigma$ on the coeficients of $ f(X)$.\\

             	       Theorem 3.7.  Assume that under the action of $Aut_k(K)$ over $(K-k)$, there is only one orbit. Then for any polynomial $f(X)$ of degree 2 in $K[X]$, $ f(K) = K.$\\

            	       Proof. Let $ f(X) = a X^2 + b X + c $. Note that $f(K) = K$ if and only if $ a^{-1}f(K) = K$. Thus we can assume that $f(X)$ is monic, hence let $ f(X) = X^2 + b X + c.$ Next note $f(K) = K $ if and only if $ (X^2 + b X)(K) = K.$ We shall show that indeed $(X^2 + b X)(K) = K.$ If $b\in k$, then the assertion follows by Theorem 3.5.  Further, by the Theorem 3.5, $ (X^2 + X + 1)(K) = K$. Thus since $ (X + 
             	       \lambda)(K) = K$ for any $ \lambda \in (K- k)$, we get that for $ g(X) = f(X)o(X+\lambda), g(K) = K.$ We have $ g(X) = X^2 + (2\lambda + 1)X + \lambda^2 + \lambda + 1$. Now, note that $ g(K) = K $ if and only if $( X^2 + (2\lambda + 1)X)(K) = K.$ Put $ d =  (2\lambda + 1) \in (K - k)$. Then $(X^2 + d X) (K) = K$. By assumption $ (K - k)/Aut_k(K) = 1$. Hence there exists $ \sigma \in Aut_k(K)$ such that $ \sigma (d) = b$. Therefore $(X^2 + d X)^{\sigma} = X^2 + b X $. Hence, by  the Remark 3.6(v), $(X^2 + b X)(K) = K$  and the result follows. \\ 
             	       
             	       Remark 3.8. From the above proof it follows that for any $ n \neq k, $ if  $ f(X)= a X^n + b X^k + c \in K[X]$, then $ f(K) = K.$\\

             	       Theorem 3.9. Let $Aut_k(K)$ has one orbit over $(K-k)$. If $ x \in (K-k)$, and $ E =  \{\alpha \in K \mid \alpha: \mbox{   algebraic over } k(x)\}$, the algebraic closure of $k(x)$ in $K$, then $ \mid (E-k)/Aut_k(E) \mid = 1$.\\
             	       
             	       Proof. Let $ y\in (E-k)$, then $y$ is transcendental over $k$ since $k$ is algebraically closed. As $ \mid (K-k)/Aut_k(K) \mid = 1$, there exists $ \sigma \in Aut_k(K)$ such that $ \sigma(x) = y$. Therefore $ \sigma (k(x)) = k(y)$. If $ \sigma(E) = F$, then $F$ is the algebraic closure of $ k(y)$  in $K$. Note that the transcendental degree of $E$ over $k$ is 1. Hence $F$ has transcendental degree 1 over $k$. As $ k \subset k(y) \subset E$ and $E$ has  transcendental degree 1 over $k$, $ E \mid k(y)$ is algebraic. Therefore $ E \subset F$. Now, as $ E\subset F, x \in (F - k)$, we can prove as above that $F\subset E$. Cosequently $ E = F$. Hence $ \sigma \in Aut_k(E)$ and $ \mid (E-k)/Aut_k(E) \mid = 1$. Thus the assertion is proved.\\
             	       
             	       We now ask the following: \\
             	       
             	       Question. Let $k$ be an algebraically closed field and $X$ an indeterminate over $k$. Let $E$ be an intermediate field such that $ k \subsetneqq E \subset \overline{k(X)}$ and  $ \mid (E-k)/Aut_k(E) \mid = 1$. Then is $ E = \overline { k(X)}$?\\
             	       
             	       Remark 3.10. In the above question, we can assume that $k(X) \subset E$, since if $y \in (E-k)$, it is transcendental over $k$. Therefore $$ k \subsetneqq k(y) \subset E \subset \overline{k(y)} = \overline{k(X)}$$ since transcedence degree of  $ \overline{k(X)}$ over $k$ is 1 ( Here $\overline{k(y)})$ is algebraic closure of $k(y)$ in $ \overline{k(X)}$). Hence we can assume $ k \subset k(X) \subset E \subset  \overline{k(X)}$ in the above question.\\
             	       
             	       Lemma 3.11. If $Aut_k(K)$ acts transitively over $(K-k)$, then for any subgroup $H$ of finite index in $ G = Aut_k(K)$,$ K^H =k$.\\
             	       
             	       Proof. If $H$ is a subgroup of finite index in $G$, then there exists a normal subgroup $H_1$ of finite index in $G$ contained in $H$. As $ K^H \subset K^{H_1}$, to prove the result we can assume $H$ is normal in $G$. Now, let $g \in G, \mbox { and }h\in H$. Then for any $ a\in K^H$,  $$ h(g(a))=  g(g^{-1}h(g(a))= g(a)$$ since $ g^{-1}h g \in H$. Cosequently $ g(K^H) \subset K^H$ for all $g\in G$. Now, as $G$ acts transitively over $(K-k), K^G = k$. Thus if $ K^H \neq k$, then $ K = K^H$. This inplies $ H $ is identity subgroup. However $G$ is infinite, hence $H$ is infinite. Therefore $ K^H = k$.\\
             	       
             	       lemma 3.12. If $Aut_k(K)$ acts transitively over $(K-k)$, then for an algebraic closure $ \overline{K} $ of $K$ either  $[\overline{K}: K] = \infty$ or $ \overline{K} = K$.\\
             	       
             	       Proof.  Assume $ 1 < [\overline{K} : K] < \infty$. Then by Artin-Screir theorem, characteristic of $K$ is $0$, and $\overline{K}  = K(i)$ where $ i^2 = -1$. By Theorem 2.4, $K$ is radically closed. Hence $ i\in K$, and $\overline{K} = K$. Thus result follows.\\  
             	       
             	       Lemma 3.13.  Assume $Aut_k(K)$ acts transitively over $(K-k)$. Let  $\overline{K}$ be an algebraic closure of $K$. Then  $\overline{K} = K$ if and only if $K$ is invariant under the action of $Aut_k(\overline{K})$.\\
             	       
             	       Proof. It is straight forward.

             	       \begin{center}   \bf{REFERENCES}
             	    \end{center}]
             	    
             	    1. Kiran S. Kedlya and Bjorn Poonen: Orbits of automorphism groups of fields, Jr. of Algebra, 293(2005), 167-184.\\
             	    
             	    2. Pramod k. Sharma : Orbits of automorphisms of integral domains, Illinois Jr. of Mathematics, Vol.52, No.2,(2008), 645-652.

                                    \end{document}